\documentclass{article}
\usepackage{amsmath,amsthm,amsfonts}
\usepackage[all]{xypic}
\newcommand{\U}{ i}
\newcommand{\ud}{{\rm d}}
\newcommand{\id}{{\rm id}}
\newcommand{\setZ}{\mathbb Z}
\newcommand{\setR}{\mathbb R}
\newcommand{\setC}{\mathbb C}
\newcommand{\g}{\mathfrak g}
\newcommand{\F}{{\mathcal F}}

\theoremstyle{plain}

\theoremstyle{definition}

\begin{document}

\title{An introduction to the Batalin-Vilkovisky formalism}
\author{Domenico Fiorenza}
\maketitle

\abstract{The aim of these notes is to introduce the quantum
master equation $\{S,S\}-2\U\hbar\Delta S=0$, and to
show its relations to the theory of Lie algebras
representations and to perturbative expansions of Gaussian
integrals. The relations of the classical master equation $\{S,S\}=0$ with
the BRST formalisms are also described.
Being an introduction, only finite-dimensional
examples will be considered.}

\section{Introduction}
The Batalin-Vilkovisky formalism is an algebraic/geometric setting to
deal with asymptotic expansions of Gaussian integrals of the form 
\[
\int_V \Psi e^{\frac{i}{\hbar}S} dv
\]
in presence of a Lie algebra $\g$ of infinitesimal symmetries of the
action $S$. From a geometric point of view, the Batalin-Vilkovisky
formalism is the theory of smooth functions on odd symplectic
supermanifolds; from an algebraic point of view is a way of looking to
Lie algebra representations (up to homotopy). 
\par
These notes are a written version of the lectures the author gave during
the seventh edition of the ``Rencontres Math\'ematiques de Glanon'', held
in
Glanon (France) in Summer 2003. They are by no means complete. Rather,
they have
to be intended as an invitation to the classical papers \cite{bat-vilk},
\cite{schwarz} and \cite{aksz}.
\par
\noindent{\it Acknowledgments} I wish to thank all the organizers and the
participants to the ``Rencontres Math\'ematiques de Glanon, \'Edition 
2003'' for the
enthusiastic and fruitful atmosphere of the meeting. I also thank Fabio
Gavarini and Riccardo Longoni who kindly read a preliminary draft of these
notes, and Pavol Severa and Jim Stasheff for their interest and useful
comments on the
first version of these notes appeared on arXiv. Special thanks go to the citizens of
Glanon, for their beautiful hospitality.

 \section{A familiar example}
The idea at the heart of the Batalin-Vilkovisky formalism is essentially
contained in the following example, which everyone is familiar with:
compute
\begin{equation*}
\int_{-\infty}^{+\infty}\frac{1}{x^2+1}\ud x
\end{equation*}
by the residues theorem. If we look closely to this computation, a lot of
interesting features appear. To begin with, to compute a 1-real variable
integral, we double the number of variables. Moreover, going from $\setR$
to
$\setC\simeq \setR^2$ something highly non-trivial has happened: the
1-form $\ud z/(z^2+1)$ extends in a non-trivial way the 1-form $\ud
x/(x^2+1)$ to a closed 1-form on ${\mathbb R}^2$. From a real variable
point of view we have changed an integral of the form
\begin{equation*}
\int_{-\infty}^{+\infty}\varphi(x)\ud x
\end{equation*}
to one of the form
\begin{equation*}
\int_{\setR\times\{0\}}\omega
\end{equation*}
where $\omega$ is a 1-form on $\setR^2$ such that
\begin{equation*}
\begin{cases}
\omega\bigr\rvert_{\setR\times\{0\}}=\varphi(x)\ud x\\
\ud \omega=0
\end{cases}
\end{equation*}
The domain of $\omega=\ud z/(z^2+1)$ is
$\setC\setminus\{\pm\U\}$. Moreover $\omega$ extends to
$\infty$, so it is actually a 1-form on ${\mathbb
P}^1(\setC)\setminus\{\pm\U\}$. The integration
domain $\setR\cup\{\infty\}$ is a cycle in ${\mathbb
P}^1(\setC)\setminus\{\pm\U\}$. Since $\omega$ is closed, we can compute
our integral by choosing another cycle in the same homology class as
$\setR\cup\{\infty\}$. If we denote by $\Gamma$ a little circle around
$\U$, then
\begin{equation*}
\int_{-\infty}^{+\infty}\frac{1}{x^2+1}\ud x
=\int_\Gamma\omega.
\end{equation*}
Now, in a little neighborhood of the point $\U$, we can expand $\omega$
into its Laurent series:
\begin{align*}
\frac{\ud
z}{z^2+1}&=\left(-\sum_{n=-1}^\infty\left(
\frac{\U}{2}\right)^{n+2}\zeta^{n}\right)\ud\zeta;\qquad\qquad
\zeta=z-\U
\end{align*}
So, if we denote by $\Gamma_0$ the circle $\Gamma-\U$, we find
\begin{align*}
\int_{-\infty}^{+\infty}\frac{1}{x^2+1}\ud x
=-\sum_{n=-1}^\infty \left(
\frac{\U}{2}\right)^{n+2}\int_{\Gamma_0}\zeta^{n}\ud\zeta
=\frac{1}{2}\int_{0}^{2\pi} \ud t
=\pi
\end{align*}

Let us analyze the various steps in the computations above.

\noindent{\em 1) The problem.} We have to compute an integral
$\int_M \Phi$,
where $M$ is an $n$-dimensional differential manifold, and $\Phi$ is a
top dimensional form on $M$.
\par
\noindent{\em 2) Doubling the coordinates.} We embed $M$ as a cycle into a
$2n$-dimensional
manifold $N$ and extend $\Phi$ to a \emph{closed} $n$-form $\Omega$ on
$N$. The condition $\ud\Omega=0$ is a \emph{gauge invariance} condition;
the manifold $M$ inside $N$ is a \emph{gauge fixing}.
\par
\noindent
{\em 3) Varying the cycle.} We choose another cycle $M_0$ in $N$, in the
same homology class as $M$. Since $\Omega$ is closed we can compute the
original integral by integrating $\Omega$ over $M_0$. Changing the
integration cycle from $M$ to $M_0$ is a \emph{change of gauge}. The fact
that the integral does not change is precisely the independence of the
result from the gauge chosen.
\par \noindent
{\em 4) The perturbative expansion.}
The cycle $M_0$ is chosen is such a way that the $n$-form $\Omega$ has a
power series expansion in a neighborhood of $M_0$. By exchanging series
and integral and obtain a series expansion for the original integral.
\par

In the example above, the closed $1$-form $\omega$ has been found by a
formal change of coordinates $x\mapsto z=x+\U\, y$. In the BV formalism,
this way of producing closed forms is called the \emph{superfield}
formalism.

\section{Super vector spaces}

A super-vector space (superspace for short) is a
$(\setZ/2\setZ)$-graded vector space $V=V_0\oplus V_1$. Here,
$V_0$ denotes the vector space of \emph{even} vectors and
$V_1$ the vector space of \emph{odd} vectors. The
$(\setZ/2\setZ)$-grading is called
\emph{degree} or \emph{parity}; the degree of a homogeneous
element $x$ will be denoted by the symbol $\overline{x}$.\par
Equivalently, one can think of a superspace as an ordinary
vector space
$V$ endowed with an automorphism $\alpha\colon V\to V$ such
that $\alpha^2=\id$. In this case, $V_0$ is the
$1$-eigenspace and $V_1$ is the $(-1)$-eigenspace. Note that
$\alpha(x)=(-1)^{\overline{x}}x$ for a homogeneous
element.\par The \emph{dimension} of a super vector space is
the pair $(\dim V_0,\dim V_1)$. If $\dim V_i=m_i$, one says
that $V$ is a $m_0|m_1$-dimensional super vector space.
All vector spaces considered in this note will be
finite-dimensional.
\par
Any vector space $V$ can be considered a super vector space
by taking $\alpha=\id$. One usually refers to this as
``placing $V$ in even degree and writes $V\oplus0$ or simply
$V$ to denote the superspace $(V,\id)$. Similarly, one can
place
$V$ in odd degree; one writes $0\oplus V[1]$ or simply
$V[1]$ to denote the superspace $(V,-\id)$.
\par
Another classical example of superspace is the following.
Given a $\setZ$-graded vector space $V=\bigoplus_{n\in
\setZ}V_n$, one can look at $V$ as the superspace with
\[
V_{\bf 0}=\bigoplus_{n\in
2\setZ}V_n;\qquad V_{\bf 1}=\bigoplus_{n\in
2\setZ+1}V_n
\]
For instance, if $(V_*,\partial)$ is a complex with a degree 1
differential, then $\partial$ can be seen as
\begin{align*}
\partial\colon V_{\bf 0}&\to V_{\bf 1}\\
\partial\colon V_{\bf 1}&\to V_{\bf 0}
\end{align*}

Now that we have defined the objects, we have to define
morphisms, in order to make superspaces a category. As one
could imagine, a morphism $\varphi\colon V\to W$ between two
superspaces is a linear map preserving the grading.
Equivalently, $\varphi$ is a linear map which intertwines
$\alpha_V^{}$ and $\alpha_W^{}$, i.e, such that the diagram
\[
\xy
\xymatrix{
   V & V \\
   W & W \\
  \ar"1,1";"1,2"^{\alpha^{}_V}
\ar"2,1";"2,2"^{\alpha^{}_W}
\ar"1,1";"2,1"_{\varphi}
\ar"1,2";"2,2"^{\varphi}
}
\endxy
\]
is commutative. \par
The category of super vector spaces (on a fixed field
${\mathbb K}$, which in these notes will always be $\setR$ or
$\setC$) will be denoted by the symbol {\sf SuperVect}.
Changing the parity of vectors in a superspace defines an
endofunctor on {\sf Supervect} which will be denoted by
$\Pi$:
\[
\Pi\colon(V,\alpha)\mapsto (V,-\alpha)
\]
that is,
\[
(\Pi V)_0=V_1;\qquad (\Pi V)_1=V_0.  
\]
We have already seen that, if $V$ is a vector space, then
$V$ can be seen as a super vector space concentrated in even
degree. The super vector space $V[1]$ obtained from $V$
concentrating it in odd degree can be seen as $\Pi V$.

 The category {\sf SuperVect} has a natural
symmetric tensor category structure defined by
\begin{align*}
(V\otimes W)_0&=(V_0\otimes W_0)\oplus (V_1\otimes W_1)\\
(V\otimes W)_1&=(V_0\otimes W_1)\oplus (V_1\otimes W_0)
\end{align*}
Equivalently,
\[
\alpha_{V\otimes W}^{}=\alpha_{V}^{}\otimes \alpha_{W}^{}
\]
The \emph{braiding} operator 
\begin{equation*}
\sigma^{}_{V,W}\colon V\otimes W\to W\otimes V
\end{equation*}
is defined as
\[
\sigma\colon x\otimes y\mapsto
(-1)^{\overline{x}\cdot\overline{y}}(y\otimes x)
\]
on homogeneous vectors $x$, $y$ and then extended by linearity.
Note that there is a natural embedding of symmetric tensor
categories
\begin{align*}
{\sf Vect}&\to{\sf Supervect}\\
V&\mapsto (V,\id)
\end{align*}
The symmetric tensor category structure on {\sf Supervect}
implies that, for any super vector space $V$ and any
positive integer $n$ there is a natural action of the
symmetric group $S_n$ on the $n^{\text{th}}$ tensor power of
$V$. The $n^{\text{th}}$ \emph{symmetric power} of $V$ is
defined as
\[
S^n(V)=(V^{\otimes n})_{S_n}=V^{\otimes
n}/\{x_1\otimes\cdots \otimes x_n-\sigma(x_1\otimes\cdots
\otimes x_n),\, \sigma \in S_n\}
\]
Note that
\[
S^n(V_0\oplus 0)=S^n(V_0)
\]
and
\[
S^n(0\oplus V_1)=\bigwedge^n V_1,
\]
so that symmetric powers of superspaces are a unified
language for symmetric and exterior powers of ordinary
vector spaces. For a general super vector space $V=V_0\oplus
V_1$ one has
\[
S^n(V_0\oplus V_1)=\bigoplus_{k=0}^n\left( S^k(V_0)\otimes
\bigwedge^{n-k}V_1\right)
\]
\section{The space of functions on a superspace}

In what follows we will be mostly concerned with the space
of regular functions on a superspace $V$. Since we are not
interested in topological questions here, we will define
regular functions as formal power series.\par
As a preliminary remark, note that the dual of a
superspace is a superspace, via
\[
\alpha_{V^*}^{}=(\alpha^{}_V)^*
\]
This implies that the even linear functionals on $V$ are
those functionals that are zero on odd vectors, and odd
linear functionals on $V$ are those functionals that are
zero on even vectors.
\par
The space of regular functions on $V$
is
\[
\F(V)=\lim_{\substack{{\longrightarrow}\\ n}}S^n(V^*)
\]
The definition of $\F(V)$ can be made more explicit by the
use of supercoordinates. Let $V$ be a $m_0|m_1$ dimensional
super vector space, and let $\{e_1,\dots, e_{m_0}\}$ be a
basis of $V_0$ and $\{\varepsilon_1,\dots,
\varepsilon_{m_1}\}$ be a basis of $V_1$. Then the dual
basis $\{x^1,\dots,x^{m_0},\theta^1,\dots,\theta^{m_1}\}$ is
a basis of $V^*$, with $x^i$ even and $\theta^k$ odd. The
linear functionals $x^i$ are called \emph{even coordinates}
on $V$ and the functionals $\theta^k$ are called \emph{odd
coordinates}. In the symmetric power $S^2(V^*)$ we have
\begin{align*}
&x^ix^j=x^jx^i\\
&x^i\theta^k=\theta^kx^i\\
&\theta^k\theta^l=-\theta^l\theta^k
\end{align*}
that is, the $x^i$ are commuting variables and the
$\theta^k$ are anticommuting (or Grassmann) variables.
With these notations,
\[
\F(V)=\setC[[x^1,\dots,x^{m_0};\theta^1,\dots \theta^{m_1}]],
\]
where $\theta^k\theta^l=-\theta^l\theta^k$.

\section{Lie algebras}

By definition, a Lie algebra is a vector space $\g$ endowed
with a bracket
\begin{equation*}
[\,,\,]\colon \g\wedge\g\to \g
\end{equation*}
which satisfies the Jacobi identity. We can look at the Lie bracket as 
a map
\begin{equation*}
[\,,\,]\colon S^2(\Pi\g)\to \Pi\g.
\end{equation*}
Let $
q:=[\,,\,]^*\colon  \Pi\g^*\to S^2(\Pi\g^*)
$
be the dual map. The space $S^2(\Pi\g^*)$ embeds into $\F(\Pi\g)$, so
we can think of $q$ as of a map
\begin{equation*}
q\colon \Pi\g^*\to \F(\Pi\g)
\end{equation*}
By forcing the Leibniz rule
$\delta(\varphi_1\varphi_2)=\delta(\varphi_1)\varphi_2
+(-1)^{\overline{\varphi_1}}\varphi_1\delta(\varphi_2)$, we
can extend
$q$ to a degree 1 derivative
\[
\delta\colon \F(\Pi\g)\to\F(\Pi\g).
\]
 The operator $\delta$ is a differential, i.e.,
 $\delta^2=0$. To see this,
we only need to show that
$\delta^2(\varphi_1\cdots\varphi_n)=0$ for any
$\varphi_1,\dots,\varphi_n\in\Pi\g^*$. Since
$\delta$ is a degree 1 derivative,
\begin{align*}
\delta(\varphi_1\cdots\varphi_n)&=
(\delta\varphi_1)\varphi_2\cdots
\varphi_n-\varphi_1(\delta\varphi_2)
\cdots\varphi_n+\\
&\qquad+\cdots+(-1)^{n-1}\varphi_1
\varphi_2\cdots
(\delta\varphi_n)
\end{align*}
Therefore
\begin{align*}
\delta^2(\varphi_1\cdots
\varphi_n)&=
(\delta^2\varphi_1)\varphi_2\cdots
\varphi_n-(\delta\varphi_1)(\delta\varphi_2)
\cdots
\varphi_n\\
&\qquad+\cdots+(-1)^{n-1}(\delta
\varphi_1)
\varphi_2\cdots(\delta\varphi_n)\\
&\qquad+(\delta\varphi_1)(\delta\varphi_2)
\cdots
\varphi_n+\varphi_1(\delta^2\varphi_2)\cdots\varphi_n\\
&\qquad+\cdots+\varphi_1\varphi_2\cdots(\delta\varphi_n)\\
&=(\delta^2\varphi_1)\varphi_2\cdots
\varphi_n+\varphi_1(\delta^2\varphi_2)\cdots
\varphi_n\\
&\qquad+\cdots+\varphi_1
\varphi_2\cdots(\delta^2\varphi_n)
\end{align*}
so, in order to prove $\delta^2=0$ we just need to prove
$\delta^2\varphi=0$ for any $\varphi\in\Pi\g^*$. By
definition,
$\delta\bigr\rvert_{\Pi\g^*}$ is the dual of the Lie bracket, i.e.,
\begin{equation*}
\langle \delta\varphi|g_1\wedge g_2\rangle=
\langle\varphi|[g_1,g_2]\rangle.
\end{equation*}
One immediately computes
\[
\langle \delta^2\varphi|g_1\wedge g_2\wedge
g_3\rangle=\langle
\varphi|[[g_1,g_2],g_3]+[g_2,g_3],g_1]+[[g_3,g_1],g_2]
\rangle
\]
So the $\delta^2=0$ is equivalent to the Jacobi identity.
\par We end this section by writing the coordinate
expression of the differential $\delta$. Let $\gamma_i$ be a
basis of $\g$ and let $c^i$ be the corresponding
coordinates on $\Pi\g$. Then
\[
(\delta c^i)(\gamma_{j_0}\wedge
\gamma_{k_0})=\langle c^i|[\gamma_{j_0},\gamma_{k_0}]\rangle=
f^i_{j_0k_0}=
\frac{1}{2}\langle f^i_{jk}c^jc^k|\gamma_{j_0}\wedge
\gamma_{k_0}\rangle,
\]
where the $f^i_{jk}$ are the structure constants of the Lie
algebra $\g$. Therefore
\[
\delta c^i=\frac{1}{2}f^i_{jk}c^jc^k
\]
that is,
\[
\delta=\frac{1}{2}f^i_{jk}c^jc^k\frac{\partial}{\partial c^i}
\]

\section{A digression on $L_\infty$-algebras}
By the discussion in the above section, we can restate the
definition of Lie algebra as follows: a Lie algebra is a
vector space
$\g$ together with a degree 1 derivative $\delta$ on
$\F(\Pi\g)$ which is a differential. if one drops the
hypothesis $\deg\delta=1$, one obtains the definition of
$L_\infty$-algebra. The structure of $L_\infty$-algebra can
also be defined by means of multilinear operations. Indeed,
if
\[
\delta\colon\F(\Pi\g)\to\F(\Pi\g)
\]
is a derivation, then $\delta$ is completely determined by
its restriction
\[
\delta\colon\Pi\g^*\to\F(\Pi\g)
\]
Let $\delta_n$ be the projection of $\delta\vert_{\Pi\g^*}$
on $S^n(\Pi\g^*)$, and let
\[
[\,,\dots,\,]_n\colon S^n(\Pi\g)\to\Pi\g
\]
be the dual maps. Then the condition $\delta^2=0$ is
equivalent to a family of quadratic relations among the
brackets $[\,,\dots,\,]_n$, and an $L_\infty$-algebra can
therefore be defined as a vector space $\g$ endowed with a
family of multilinear brackets
\[
[\,,\dots,\,]^{}_n\colon\bigwedge^n\g\to\g
\]
with $\deg[\,,\dots,\,]^{}_n=n-1$ and satisfying certain
quadratic relations. It is interesting to write down the
first of these relations. With the above notations,
$\delta^2=(\delta_1+\delta_2+\cdots)^2=\delta_1^2+
(\delta_1\delta_2+\delta_2\delta_1)+\cdots$, so that
$\delta^2=0$ implies
\begin{align*}
\delta_1^2&=0\\
\delta_1\delta_2&+\delta_2\delta_1=0\\
\delta_1\delta_3+&\delta_2^2+\delta_3\delta_1=0
\end{align*}
and so on. Consider now the brackets
\begin{align*}
[\,]^{}_1=\delta_1^*\colon\g&\to\g\\
[\,]^{}_2=\delta_2^*\colon\g\wedge\g&\to\g\\
[\,]^{}_3=\delta_3^*\colon\g\wedge\g\wedge\g&\to\g\\
\end{align*}
The above equations for the $\delta_i$ tell us that $[\,]^{}_1$
is a differential on $\g$. Moreover, if we set $d=[\,]^{}_1$,
then
\[
d[g_1,g_2]^{}_2=[dg_1,g_2]^{}_2+[g_1,dg_2]^{}_2
\]
i.e., $d$ is a derivative with respect to the bracket
$[\,]^{}_2$. Finally,
\begin{multline*}
[[g_1,g_2]_2,g_3]^{}_2+
[[g_2,g_3]_2,g_1]^{}_2+
[[g_3,g_1]_2,g_2]^{}_2=\\d[g_1,g_2,g_3]^{}_3-[dg_1,g_2,g_3]^{}_3-
[g_1,dg_2,g_3]^{}_3-[g_1,dg_2,g_3]^{}_3
\end{multline*}
i.e., the Jacoby relation for the bracket $[\,]^{}_2$ holds up
to a homotopy given by the bracket $[\,]^{}_3$. In this sense
an $L_\infty$-algebra is a Lie algebra up to homotopy.
Moreover, since the homotopies $[\,]^{}_n$ with $n\geq 3$ are
part of the data defining the structure of $L_\infty$-algebra,
$L_\infty$-algebras are
also called \emph{strong homotopy Lie algebras}.
\section{Lie algebras representations}
Let $V$ be a vector space and $\g$ a Lie algebra. By
definition, a representation of $\g$ on $V$ is a linear map
\[
\rho\colon \g\otimes V\to V
\]
such that the induced map
\[
\g\to\text{End}(V)
\]
is a Lie algebra morphism. Equivalently, the diagram
\[
\xy
\xymatrix{
   (\g\wedge\g)\otimes V & \g\otimes(\g\otimes
V) &\g\otimes V
\\
   \g\otimes V && V \\
  \ar"1,1";"1,2"^{\iota}
\ar"2,1";"2,3"^{\rho}
\ar"1,1";"2,1"_{[\,,\,]}
\ar"1,3";"2,3"^{\rho}
\ar"1,2";"1,3"^{\phantom{mmi}\id\otimes\rho}
}
\endxy
\]
is commutative, where the inclusion  $\iota\colon\g\wedge\g
\hookrightarrow\g\otimes \g$ is given by $g_1\wedge
g_2\mapsto g_1\otimes g_2-g_2\otimes g_1$.
The three maps
\begin{align*}
0\colon S^2(V)&\to 0\\
\rho\colon \g \otimes V&\to V\\
[\,,\,]\colon \g\wedge\g&\to \g
\end{align*}
can be seen as a map
\[
\rho+[\,,\,]\colon S^2(V\oplus\Pi\g)\to V\oplus\Pi\g
\]
by the isomorphism
\[
S^2(V\oplus\Pi\g)\simeq S^2(V)\oplus (\Pi\g\otimes V)\oplus
S^2(\Pi\g)
\]
Denote by $\delta$ the map dual to $\rho+[\,,\,]$,
\[
\delta\colon
(V\oplus\Pi\g)^*\to S^2\bigl((V\oplus\Pi\g)^*\bigr)
\]
and extend it to a degree 1 derivation
\[
\delta\colon
\F(V\oplus\Pi\g)\to\F(V\oplus\Pi\g)
\]
by forcing the Leibniz rule. Reasoning as for Lie algebras,
one sees that $\delta$ is a differential, i.e.,
$\delta^2=0$. Indeed, for any $\varphi\in (V\oplus\Pi\g)^*$,
\begin{align*}
\langle
\delta^2\varphi|v_1v_2v_3\rangle&=\langle\varphi|0\rangle
\\
\langle \delta^2\varphi|g\otimes
v_1v_2\rangle&=\langle\varphi|0\rangle
\\
\langle \delta^2\varphi|(g_1\wedge g_2)\otimes
v\rangle&=\langle\varphi|[g_1,g_2]\cdot v
-g_1\cdot(g_2\cdot v)+g_2\cdot(g_1\cdot v)\rangle
\\
\langle \delta^2\varphi|g_1\wedge g_2\wedge
g_3\rangle&=\langle\varphi|[[g_1,g_2],g_3]+[g_2,g_3],g_1]
+[[g_3,g_1],g_2]\rangle
\end{align*}
where we have written $g\cdot v$ for $\rho(g\otimes v)$.
These equations show that the condition $\delta^2=0$ is
equivalent to the Jacobi identity for $\g$ and to the fact
that $\rho$ is a representation. That is, a Lie algebra
representation $(\g,V,\rho,[\,,\,])$ can be seen as a
superspace $V\oplus\Pi\g$ endowed with a degree one
derivative $\delta\colon \F(V\oplus\Pi\g)\to
\F(V\oplus\Pi\g)$ which is a differential, i.e.,
$\delta^2=0$.\par
The operator $\delta$ is called the \emph{BRST operator},
and the cohomology of $\delta$ is called the \emph{BRST
cohomology}. Here BRST is the acronym for Becchi-Rouet-Stora-Tyutin.

A more refined analysis of BRST cohomology can be obtained
taking care of the degrees. To begin with, recall that
\[
\F(V\oplus\Pi\g)=\bigoplus_{p,q}S^p(V^*)\otimes
S^q(\Pi\g^*)\simeq\bigoplus_{p,q}S^p(V^*)\otimes
\bigwedge^q(\g^*)
\]
so that $\delta$ can be seen as an operator
\[
\delta^{p,q}\colon S^p(V^*)\otimes \bigwedge^q\g^*\to
S^p(V^*)\otimes \bigwedge^{q+1}\g^*
\]
In particular, if we take $p=0$ we obtain
\[
\delta^{0,q}\colon \bigwedge^q\g^*\to
\bigwedge^{q+1}\g^*
\]
Dualizing this differential we obtain a differential
\[
d_{q+1}\colon \bigwedge^{q+1}\g\to
\bigwedge^{q}\g
\]
which is easily seen to be the Chevalley-Eilenberg
differential defining the Lie algebra cohomology of $\g$.
Another classical example is obtained by taking $p=1$. In
this case one has
\[
\delta^{1,q}\colon V^*\otimes\bigwedge^q\g^*\to
V^*\otimes\bigwedge^{q+1}\g^*
\]
Dualizing this differential we obtain a differential
\[
d_{q+1}\colon V\otimes\bigwedge^{q+1}\g\to
V\otimes \bigwedge^{q}\g
\]
which is the
differential defining the Lie algebra cohomology of $\g$
with coefficients in the representation $\rho$.
\par As we did for Lie algebras, we end this section by
writing the coordinate expression of the differential
$\delta$. Let $e_i$ and
$\gamma_j$ be  basis of $V$ and $\g$ respectively and let
$v^i$ and
$c^i$ be the corresponding coordinates on $V$ and $\Pi\g$.
Then
\begin{align*}
(\delta v^i)(e_{j_0}\cdot
e_{k_0}&+e_{j_1}\otimes \gamma_{k_1}+\gamma_{j_2}\wedge
\gamma_{k_2})
=\langle
v^i|\rho(\gamma_{k_1})e_{j_1}+[\gamma_{j_2},\gamma_{k_2}]
\rangle\\
&= \rho^i_{j_1k_1}=
\langle \rho^i_{jk}v^jc^k|e_{j_0}\cdot
e_{k_0}+e_{j_1}\otimes \gamma_{k_1}+\gamma_{j_2}\wedge
\gamma_{k_2}\rangle,
\end{align*}
where the $\rho^i_{jk}$ are the structure constants of the Lie
algebra representation $\rho\colon\g\to\text{End}(V)$.
The computation we already did for Lie algebras give
\begin{align*}
(\delta c^i)(e_{j_0}\cdot
e_{k_0}&+e_{j_1}\otimes \gamma_{k_1}+\gamma_{j_2}\wedge
\gamma_{k_2})=\langle
c^i|\rho(\gamma_{k_1})e_{j_1}+[\gamma_{j_2},\gamma_{k_2}]
\rangle\\
&=
f^i_{j_2k_2}=
\frac{1}{2}\langle f^i_{jk}c^jc^k|e_{j_0}\cdot
e_{k_0}+e_{j_1}\otimes \gamma_{k_1}+\gamma_{j_2}\wedge
\gamma_{k_2}\rangle,
\end{align*}
where the $f^i_{jk}$ are the structure constants of the Lie
algebra $\g$. Therefore
\[
\delta v^i=\rho^i_{jk}v^jc^k ;\qquad \delta
c^i=\frac{1}{2}f^i_{jk}c^jc^k
\]
that is,
\[
\delta=\rho^i_{jk}v^jc^k\frac{\partial}{\partial
v^i}+
\frac{1}{2}f^i_{jk}c^jc^k\frac{\partial}{\partial
c^i}
\]

\section{Batalin-Vilkovisky algebras}

Let now $W$ be any superspace. The superspace $W\oplus\Pi
W^*$ is naturally endowed with an odd non-degenerate
pairing. Let $\Delta$ be the Laplace operator associated to
this pairing. If $x^i$ are coordinates on $W$ (the \emph{fields}) and
$x^+_i$ are the corresponding coordinates on $\Pi W^*$ (the
\emph{antifields}), then
\[
\Delta=\frac{\partial}{\partial
x_i^+}\frac{\partial}{\partial x^i}
\]
The operator $\Delta$ is called the Batalin-Vilkovisky
Laplacian; note that, if $\Phi$ is a
homogeneous function in $\F(W\oplus\Pi W^*)$, then
$\Delta\Phi$ is also homogeneous and
$\overline{\Delta\Phi}=\overline{\Phi}+1\mod 2$. It is
immediate to compute that
\[
\Delta^2=0
\]
Indeed,
\begin{align*}
\Delta^2&=\frac{\partial}{\partial
x^+_i}\frac{\partial}{\partial x^i}\frac{\partial}{\partial
x^+_j}\frac{\partial}{\partial
x^j}\\
&=(-1)^{\overline{x^i}\cdot\overline{x^+_j}+
\overline{x^i}\cdot\overline{x^j}+
\overline{x^+_i}\cdot\overline{x^+_j}+\overline{x^+_i}\cdot
\overline{x^j}}
\frac{\partial}{\partial
x^+_j}\frac{\partial}{\partial x^j}\frac{\partial}{\partial
x^+_i}\frac{\partial}{\partial
x^i}\\
&=(-1)^{(\overline{x^i}+\overline{x^+_i})(
\overline{x^j}+\overline{x^+_j})}
\frac{\partial}{\partial
x^+_j}\frac{\partial}{\partial x^j}\frac{\partial}{\partial
x^+_i}\frac{\partial}{\partial
x^i}.
\end{align*}
Since the variables $x^i$ and $x^+_i$ have opposite parity,
$\overline{x^i}+\overline{x^+_i}=1\mod 2$, for any $i$.
Therefore,
\[
\Delta^2=-\Delta^2
\]
i.e., $\Delta^2=0$. The cohomology of $\F(W\oplus\Pi W^*)$ with respect 
to the BV-Laplacian is called \emph{BV-cohomology} or
$\Delta$-\emph{cohomology}.
\par
Let now
$\Phi$ and $\Psi$ be two homogeneous functions on $W\oplus\Pi W^*$. Then
\begin{align*}
\Delta(\Phi\cdot\Psi)&=\frac{\partial}{\partial
x^+_i}\frac{\partial}{\partial x^i}(\Phi\cdot \Psi)\\
&=\frac{\partial}{\partial
x^+_i}\left(\frac{\partial\Phi}{\partial x^i}\cdot
\Psi+(-1)^{\overline{x^i}\cdot\overline{\Phi}}
\Phi\frac{\partial\Psi}{\partial
x^i}\right)\\
&=\frac{\partial}{\partial
x^+_i}\frac{\partial\Phi}{\partial x^i}\cdot
\Psi
+(-1)^{(\overline{x^i}+\overline{\Phi})\overline{x^+_i}}
\frac{\partial\Phi}{\partial x^i}
\frac{\partial\Psi}{\partial
x^+_i}
+(-1)^{\overline{x^i}\cdot\overline{\Phi}}
\frac{\partial\Phi}{\partial
x^+_i}\frac{\partial\Psi}{\partial
x^i}+\\
&\phantom{\frac{\partial}{\partial
x^+_i}\frac{\partial\Phi}{\partial x^i}\cdot
\Psi
+(-1)^{(\overline{x^i}+\overline{\Phi})\overline{x^+_i}}
\frac{\partial\Phi}{\partial x^i}
\frac{\partial\Psi}{\partial
x^+_i}}+(-1)^{(\overline{x^i}+\overline{x^+_i})
\overline{\Phi}}\Phi
\frac{\partial}{\partial
x^+_i}\frac{\partial\Psi}{\partial
x^i}\\
&=
(\Delta\Phi)\cdot\Psi+(-1)^{\overline{\Phi}}\{\Phi,\Psi\}
+(-1)^{\overline{\Phi}}
\Phi\cdot\Delta\Psi
\end{align*}
where $\{\Phi,\Psi\}$ is the so-called BV-bracket, defined by
\[
\{\Phi,\Psi\}=(-1)^{\overline{x^+_i}\cdot\overline{\Phi}}
\frac{\partial\Phi}{\partial
x^+_i}\frac{\partial\Psi}{\partial
x^i}
-(-1)^{(\overline{\Phi}+1)(\overline{\Psi}+1)
+\overline{x^+_i}\cdot\overline{\Psi}
}
\frac{\partial\Psi}{\partial
x^+_i}
\frac{\partial\Phi}{\partial x^i}
\]
The BV-bracket is best expressed using \emph{left} and
\emph{right derivatives}: for a homogeneous vector $v$ in
$W\oplus\Pi W^*$, set
\begin{align*}
\overrightarrow{\partial}_v\Phi&=\partial_v\Phi\\
\overleftarrow{\partial}_v\Phi&=(-1)^{\overline{v}
\cdot\overline{\Phi}}\partial_v\Phi
\end{align*}
With these notations, the BV-bracket reads
\[
\{\Phi,\Psi\}=
\frac{\overleftarrow{\partial}\Phi}{\partial
x^+_i}\frac{\overrightarrow{\partial}\Psi}{\partial
x^i}
-(-1)^{(\overline{\Phi}+1)(\overline{\Psi}+1)}
\frac{\overleftarrow{\partial}\Psi}{\partial
x^+_i}
\frac{\overrightarrow{\partial}\Phi}{\partial x^i}
\]
The Leibniz rule for derivatives gives
\begin{align*}
\{\Phi,\Psi\cdot \Upsilon\}&=
(-1)^{\overline{x^+_i}\cdot\overline{\Phi}}
\frac{\partial\Phi}{\partial
x^+_i}\frac{\partial\Psi}{\partial
x^i}\Upsilon+(-1)^{\overline{x^+_i}\cdot\overline{\Phi}
+(\overline{\Phi}+1)\overline{\Psi}
}
\Psi\frac{\partial\Phi}{\partial
x^+_i}\frac{\partial\Upsilon}{\partial
x^i}+\\
&\phantom{m}
-(-1)^{(\overline{\Phi}+1)(
\overline{\Psi}+\overline{\Upsilon}+1)+\overline{x^+_i}
(\overline{\Psi}+\overline{\Upsilon})}
\frac{\partial\Psi}{\partial
x^+_i}\Upsilon\frac{\partial\Phi}{\partial
x^i}+\\
&\phantom{m}
-(-1)^{(\overline{\Phi}+1)(
\overline{\Psi}+\overline{\Upsilon}+1)+\overline{x^+_i}
(\overline{\Psi}+\overline{\Upsilon})+\overline{x^+_i}
\cdot\overline{\Upsilon}}
\Psi\frac{\partial\Upsilon}{\partial
x^+_i}\frac{\partial\Phi}{\partial
x^i}
\\
&=\{\Phi,\Psi\}\Upsilon+(-1)^{(\overline{\Phi}+1)\overline{\Psi}}
\Psi\{\Phi,\Upsilon\}
\end{align*}
i.e., the BV-bracket satisfies an odd Poisson identity: for
any homogeneous $\Phi\in\F(W\oplus\Pi W^*)$,
the operator $\text{ad}_\Phi=\{\Phi,-\}$ is a derivative
of degree
$\overline{\Phi}+1$.
\par
The data $\F(W\oplus\Pi W^*)$, $\cdot$, $\Delta$ and
$\{\,,\,\}$ are the basic example of \emph{BV-algebra}. More
generally, a BV-algebra is $({\mathcal
A},\cdot,\Delta,\{\,,\,\})$, where
\begin{enumerate}
\item ${\mathcal A}$ is a superspace;
\item $\cdot\colon {\mathcal A}\otimes {\mathcal A}\to
{\mathcal A}$ is an associative and graded-commutative
multiplication;
\item $\Delta\colon {\mathcal A}\to {\mathcal A}$ is an odd
differential (here odd means that $\Delta$ changes the
parity of a homogeneous element);
\item $\{\,,\,\}\colon {\mathcal A}\otimes {\mathcal A}\to
{\mathcal A}$ is a bilinear operation such that
\[
\Delta(\Phi\cdot\Psi)=(\Delta\Phi)\cdot\Psi+(-1)^{\overline{\Phi}}\{\Phi,\Psi\}
+(-1)^{\overline{\Phi}}
\Phi\cdot\Delta\Psi
\]
and
\[\{\Phi,\Psi\cdot \Upsilon\}=
\{\Phi,\Psi\}\Upsilon+(-1)^{(\overline{\Phi}+1)\overline{\Psi}}
\Psi\{\Phi,\Upsilon\}
\]
\end{enumerate}
The operations $\cdot$, $\Delta$ and
$\{\,,\,\}$ are called, respectively, \emph{multiplication},
\emph{BV-Laplacian} and \emph{BV-bracket}. Note that the
parity of $\{\Phi,\Psi\}$ is
$\overline{\Phi}+\overline{\Psi}+1\mod 2$.

We will now derive some consequences of the compatibility
between the three operations in a BV-algebra. To begin
with, the BV-bracket is anticommutative, up to a parity
change:
\[
\{\Psi,\Phi\}=-(-1)^{(\overline{\Psi}+1)(\overline{\Phi}+1)}
\{\Phi,\Psi\}
\]
Indeed,
$\Psi\cdot\Phi=(-1)^{\overline{\Psi}\cdot\overline{\Phi}}
\Phi\cdot\Psi$, so that
\begin{align*}
\Delta(\Psi\cdot\Phi)&=(-1)^{\overline{\Psi}
\cdot\overline{\Phi}}
\Delta(\Phi\cdot\Psi)\\
&=
(-1)^{\overline{\Psi}
\cdot\overline{\Phi}}(\Delta\Phi)\cdot\Psi+
(-1)^{\overline{\Phi}(\overline{\Psi}+1)}\{\Phi,\Psi\}
+(-1)^{\overline{\Phi}(\overline{\Psi}+1)}
\Phi\cdot\Delta\Psi
\end{align*}
Since the BV-Laplacian changes the parity of a homogeneous
function, we can rewrite this as
\begin{align*}
\Delta(\Psi\cdot\Phi)&=
(-1)^{\overline{\Psi}
}\Psi\cdot \Delta\Phi+
(-1)^{\overline{\Phi}(\overline{\Psi}+1)}\{\Phi,\Psi\}
+(\Delta\Psi)\cdot\Phi
\end{align*}
On the other hand,
\[
\Delta(\Psi\cdot\Phi)=
(\Delta\Psi)\cdot\Phi+(-1)^{\overline{\Psi}}\{\Psi,\Phi\}
+(-1)^{\overline{\Psi}}
\Psi\cdot\Delta\Phi
\]
so that, equating the right hand sides of the two equations
above, we find
\[
(-1)^{\overline{\Phi}(\overline{\Psi}+1)}\{\Phi,\Psi\}=
(-1)^{\overline{\Psi}}\{\Psi,\Phi\}
\]
which is the equation we were looking for. Note that odd
anticommutativity for the BV-bracket on $\F(W\oplus\Pi W^*)$
could be immediately checked from the definition of the
bracket.
\par
The compatibility equation between the three operations in a
BV-algebra shows that the Laplacian is not a derivative with
respect to the multiplication: indeed, the BV-bracket measures
precisely how much $\Delta$ fails to be a derivative. Yet,
the Laplacian is a derivative with respect to the bracket:
\[
\Delta\{\Phi,\Psi\}=\{\Delta\Phi,\Psi\}+
(-1)^{\overline{\Phi}+1}\{\Phi,\Delta\Psi\}
\]
This is an immediate consequence of the trivial identity
$\Delta^2(\Phi\cdot\Psi)=0$. Finally, applying the
BV-Laplacian to the odd Poisson identity $\{\Phi,\Psi\cdot \Upsilon\}=
\{\Phi,\Psi\}\Upsilon+(-1)^{(\overline{\Phi}+1)\overline{\Psi}}
\Psi\{\Phi,\Upsilon\}$ one finds the odd Jacobi identity for the
BV-bracket:
\[
\{\Phi,\{\Psi,\Upsilon\}\}=\{\{\Phi,\Psi\},\Upsilon\}+(-1)^{
(\overline{\Phi}+1)(\overline{\Psi}+1)}\{\Psi,\{\Phi,\Upsilon\}\}
\]
\par
As a concluding remark, note that by the compatibility equation between
the BV-Laplacian, the multiplication and the BV-bracket
\[
\Delta(\Phi\cdot\Psi)=(\Delta\Phi)\cdot\Psi+(-1)^{\overline{\Phi}}\{\Phi,\Psi\}
+(-1)^{\overline{\Phi}}
\Phi\cdot\Delta\Psi
\]
one can express the BV-bracket entirely in terms of the of the
multiplication and of the  BV-Laplacian. The odd Poisson identity is then
translated into the following \emph{seven-terms relation}:
\begin{multline*}
\Delta(\Phi\Psi\Upsilon)+(\Delta\Phi)\Psi\Upsilon+
(-1)^{\overline{\Phi}}\Phi(\Delta\Psi)\Upsilon+(-1)^{\overline{\Phi}+
\overline{\Psi}}\Phi\Psi(\Delta\Upsilon)=\\
\Delta(\Phi\Psi)\Upsilon+(-1)^{\overline{\Phi}}\Phi\Delta(\Psi\Upsilon)
+(-1)^{(\overline{\Phi}+1)
\overline{\Psi}}\Psi\Delta(\Phi\Upsilon)
\end{multline*}
Therefore, one could define BV algebras as the data $({\mathcal
A},\cdot,\Delta,)$, where ${\mathcal A}$ is a superspace, $\cdot\colon
{\mathcal A}\otimes {\mathcal A}\to
{\mathcal A}$ is an associative and graded-commutative
multiplication, and $\Delta\colon {\mathcal A}\to {\mathcal A}$ is an odd
differential such that the seven-terms relation hold. The BV-bracket would
then be defined by the formula
\[
\{\Phi,\Psi\}
=(-1)^{\overline{\Phi}}\Delta(\Phi\cdot\Psi)+
(-1)^{\overline{\Phi}+1}(\Delta\Phi)\cdot\Psi
-\Phi\cdot\Delta\Psi
\]
On the other hand, the BV-bracket is such an important operation in the 
Batalin-Vilkovisky formalism, that we preferred to make it part of the
definition of BV-algebra.

\section{BV cohomology, Lagrangian submanifolds and the
quantum master equation}

This is the most geometrical part of this short note. Being just an
introduction, we will be very sketchy and invite the reader to look at
the details into \cite{schwarz} and \cite{aksz}. The starting point is
that, for any superspace $W$, the superspace $W\oplus \Pi W^*$ has a
canonical odd symplectic structure. This is very familiar from
classic differential geometry: for any vector space $V$, the vector space
$V\oplus V^*$ has a canonical symplectic structure. Since we are dealing
with an (odd) symplectic space, it is meaningful to consider Lagrangian
submanifolds of $W\oplus\Pi W^*$: these are just
isotropic (super-)submanifolds of maximal dimension. 
The volume form $dx^1 dx^2\cdots dx^ndx_1^+ dx_2^+\dots
dx_n^+$ induces well defined volume forms on the Lagrangian
submanifolds of $W\oplus\Pi W^*$, so the functionals
\[
\int_{\mathcal L}
\]
 on $\F(W\oplus\Pi W^*)$ are defined for any Lagrangian ${\mathcal
L}\subseteq W\oplus \Pi W^*$. Here comes the main theorem of the
BV-formalism. Let $\Phi\in\F(W\oplus\Pi W^*)$. 
\[
\boxed{
\begin{matrix}
\text{ If $\Delta\Phi=0$, then }
\displaystyle{\int_{\mathcal L}\Phi} 
\text{ depends only on the homology class of ${\mathcal L}$.}\\
\text{Moreover, if $\Phi=\Delta\Psi$, then } \displaystyle{\int_{\mathcal
L}\Phi=0}\text{ for
any Lagrangian ${\mathcal L}$.}
\end{matrix}
}
\]
This can be thought as a BV version of the familiar Stokes'
theorem. In particular, we have a pairing between homology classes of
Lagrangian submanifolds and $\Delta$-cohomology classes of functions on
the superspace $W\oplus \Pi W^*$.

  As in classical
symplectic manifolds, for any smooth function
$F$ on $W$, the sub-manifold ${\mathcal L}_F$ defined by the equations
\[
x^+_i=\frac{\partial F}{\partial x^i}
\]
is Lagrangian. Moreover, for any two functions $F_0$ and $F_1$ the
homotopy $F_t=tF_1+(1-t)F_0$ shows that ${\mathcal L}_{F_0}$ and
${\mathcal L}_{F_1}$ are in the same homology class. Therefore
\[
\int_{{\mathcal L}_{F_0}}\Phi=\int_{{\mathcal L}_{F_1}}\Phi
\]
for any $\Delta$-closed $\Phi$. In particular, the subspace $W$ of
$W\oplus \Pi W^*$ is defined by the equations
\[
x_i^+=0
\]
so it is the Lagrangian subspace defined by the function $F\equiv 0$. It
follows that, for any $\Delta$-closed function $\Phi$ and any smooth
function $F$,
\[
\int_{W}\Phi=\int_{{\mathcal L}_{F}}\Phi.
\]
The function $F$ is called the \emph{gauge fixing fermion}.
 
Integrands we are interested in are usually of the form
\[
\Phi=\Psi e^{\frac{i}{\hbar}S}
\]
We begin by assuming that $\Psi=1$, so that the equation
$\Delta\Phi=0$ becomes
\[
0=\Delta e^{\frac{i}{\hbar}S}=\Delta\left(
\sum_{n=0}^\infty \frac{(iS)^n}{\hbar^n n!}
\right)
=\left(\frac{i}{\hbar}\Delta
S-\frac{1}{2\hbar^2}\{S,S\}\right)e^{\frac{i}{\hbar}S}
\]
so the condition $\Delta e^{\frac{i}{\hbar}S}=0$ is
equivalent to the \emph{quantum master equation}
\[
\boxed{
\{S,S\}-2i\hbar\Delta S=0
}
\]
Note that, since $S$ is even, the bracket $\{S,S\}$ is non
trivial. If $S$ can be expanded as a series in $\hbar$,
\[
S=S_0+\hbar S_1+\hbar^2 S_2+\cdots
\]
then the quantum master equation becomes the sequence of
equations
\begin{align*}
\{S_0,S_0\}&=0\\
\{S_0,S_1\}&=i\Delta S_0\\
\{S_0,S_2\}+\frac{1}{2}\{S_1,S_1\}&=i\Delta S_1\\
&\cdots
\end{align*}
The first equation in the above list is called the
\emph{master equation}; since $S_0$ is even, it is a
non-trivial equation. We will see in the next section how to
relate solutions of the master equation with Lie algebra
representations. But now, let us go back to the problem of
calculating the Gaussian integral
\[
\langle\!\langle \Psi\rangle\!\rangle=\int_W \Psi\,
e^{\frac{i}{\hbar}S}
\]
via the BV formalism.
Assume that $S$ is a solution of the quantum master equation.
In order to apply the BV formalism, we want
\[
\Delta\left(\Psi\, e^{\frac{i}{\hbar}S}\right)=0
\]
which is equivalent to
\[
\Delta\Psi+\frac{i}{\hbar}\{S,\Psi\}=0
\]
Let $\Omega$ be the operator
\[
\Omega=-i\hbar\Delta+\text{ad}_S
\]
then the BV formalism applies to all the functions $\Psi$ in
$\ker\Omega$. The operator $\Omega$ is actually a
differential; indeed,
\begin{align*}
\Omega^2\Psi&=\Omega\left(-i\hbar\Delta\Psi+\{S,\Psi\}\right)\\
&=-\hbar^2\Delta^2\Psi-i\hbar\Delta\{S,\Psi\}
-i\hbar\{S,\Delta\Psi\}+
\{S,\{S,\Psi\}\}\\
&=-i\hbar\{\Delta S,\Psi\}+\frac{1}{2}
\{\{S,S\},\Psi\}\\
&=\frac{1}{2}\{\{S,S\}-2i\hbar\Delta
S,\Psi\}=0.
\end{align*}
Moreover, if $\Psi$ is $\Omega$-exact, then
\begin{align*}
\langle\!\langle \Psi\rangle\!\rangle&=
\int_W (\Omega\Psi_0)\,
e^{\frac{i}{\hbar}S}=-i\hbar\int_W \Delta\left(\Psi_0\,
e^{\frac{i}{\hbar}S}\right)=0
\end{align*}
Therefore, the expectation value is actually a linear
functional on the $\Omega$-cohomology of $\F(W\oplus\Pi
W^*)$; the $\Omega$-cohomology classes are called
\emph{observables} of the theory. The operator $\Omega$ and
its cohomology are called, respectively, the quantum
BRST operator and the quantum BRST cohomology. 

\section{From Lie algebras representations to solutions of
the master equation}
Recall that to a Lie algebra representation $\rho\colon
\g\otimes V\to V$ is associated a BRST differential $\delta$
on $\F(V\oplus\Pi\g)$. Set
\[
W=V\oplus \Pi\g
\]
and consider the superspace
\[
W\oplus\Pi W^*
\]
As we have seen, $\F(W\oplus\Pi W^*)$ has a natural
structure of BV-algebra. Let $S_0$ be any element in
$\F(W)$. It can be considered as an element
of
$\F(W\oplus\Pi W^*)$ by the projection $W\oplus \Pi W^*\to
W$. Clearly, as an element of $\F(W\oplus\Pi W^*)$, the
function $S_0$ satisfies the master equation:
$\{S_0,S_0\}=0$; indeed,
$S_0$ is independent of the coordinates of vectors in $\Pi
W^*$. Therefore, the operator $\text{ad}_{S_0}=\{S_0,-\}$ is
a differential on $\F(W\oplus\Pi W^*)$ which induces the
zero differential on $\F(W)$. Moreover, if $\g$ is a Lie
algebra of infinitesimal symmetries of the action $S_0$,
then $\delta S_0=0$.
\par
Let now $S_1$ be the quadratic
function on $W\oplus\Pi W^*$ defined by
\[
S_1(w+\omega)=\langle\delta w|\omega\rangle
\]
If $x^1,\dots,x^n$ are coordinates on $W$ and $x_i^+,\dots,
x_n^+$ are the dual coordinates on $\Pi W^*$, then
\[
S_1=x_i^+\delta x^i
\]
Clearly,
\[
\{S_1,\Phi\}=\delta \Phi
\]
for any $\Phi\in\F(W)$. Moreover,
\[
\delta S_1=(-1)^{\overline{x}_i^+}x^+_i\delta^2 x^i=0,
\]
since $\delta^2=0$.
We compute
\begin{align*}
\{S_1,S_1\}&=\{S_1,x_i^+\delta
x^i\}\\
&=\{S_1,x_i^+\}\delta
x^i+(-1)^{\overline{x}_i^+} x_i^+\{S_1,\delta x^j\}\\
&=\frac{\partial S_1}{\partial x^i}\delta x^1
+(-1)^{\overline{x}_i^+}
x_i^+\delta^2 x^i\\
&=2\delta S_1=0
\end{align*}
Since $\{S_1,S_0\}=\delta S_0$, if $\g$ is a Lie
algebra of infinitesimal symmetries of the action $S_0$, then
$\{S_1,S_0\}=0$. So, if we set
\[
S=S_0+\hbar S_1
\]
then
\[
\{S,S\}=0
\]
that is $S$ is a solution of the master equation. As a
consequence,
$\text{ad}_S=\{S,-\}$ is a differential on
$\F(W\oplus\Pi W^*)$. Moreover $\{S,-\}$ induces the
differential $\hbar\delta$ on $\F(W)$, since
\[
\{S,\Phi\}=\{S_0,\Phi\}+\hbar\{S_1,\Phi\}=
\hbar\{S_1,\Phi\}=\hbar\delta\Phi
\]
for any $\Phi\in\F(W)$.

\section{On-shell and off-shell representations of Lie
algebras}

In the previous section we have seen how to write a solution
to the master equation starting from a Lie algebra
representation and a function $S_0$ for which the Lie
algebra $\g$ was an algebra of infinitesimal symmetries.
Moreover, the solution
$S=S_0+\hbar S_1$ was first order in the antifields. We will
now show how to go the other way round, i.e., how to
associate a Lie algebra representation and a function $S_0$
for which the Lie algebra is an algebra of infinitesimal
symmetries to a solution to the master equation which is
first order in the antifields.
\par
Let $S=S_0+\hbar S_1$ be a solution of the master equation.
We stress the fact that $S_k$ is of degree $k$ in the
antifields, for $k=0,1$ by writing
\[
S=S_0(x^1,\dots, x^n)+\hbar x_i^+Q^i(x^1,\dots,x^n)
\]
Since $S_0$ does not depend on the antifields, the master
equation $\{S,S\}$ is reduced to the two equations
\begin{align*}
\{S_1,S_0\}&=0\\
\{S_1,S_1\}&=0
\end{align*}
The second equation of this pair
tells that
$\text{ad}_{S_1}$ is a differential on $\F(W\oplus\Pi W^*)$.
Moreover, since $S_1$ is first-order in the antifields,
$\text{ad}_{S_1}$ maps $\F(W)$ to itself. Let
$\delta\colon
\F(W)\to
\F(W)$ be the differential defined by the restriction of
$\text{ad}_{S_1}$ to $\F(W)$. 
It is immediate to compute
\[
\delta x^i=\{S_1,x^i\}=Q^i(x^1,\dots,x^n)
\]
so that, if $\deg Q_i=2$, then $\delta$ is a degree one
derivative on $\F(W)$ which is a differential. As we have
shown in section 7, this defines a Lie algebra representation
$\g\to\text{End}(V)$, where $V$ and $\g$ are defined by
$W=V\oplus\Pi\g$, i.e., $V=W_0$ and $\g=\Pi W_1$. Note that
the equations $\delta x^i=Q^i(x^1,\dots,x^n)$ imply that
$S_1=x_i^+\delta x^i$, i.e., $S_1$ has the form seen in
section 10.
\par
Finally, the equation $\{S_1,S_0\}$ gives $\delta S_0=0$,
that is, the Lie algebra $\g$ is an algebra of infinitesimal
symmetries for $S_0$. If we drop the condition $\deg Q_i=2$
then we get into the realm of homotopy Lie algebras
representations, see section 6 above. 
\par
If we drop the hypothesis $S$ being of degree 1
in the antifields, and write 
\[
S=S_0+\hbar S_1+\hbar^2 S_2+\cdots
\]
with $S_k$ of degree $k$ in the antifields, then the master
equation $\{S,S\}=0$ becomes the sequence of equations
\begin{align*}
\{S_0,S_1\}&=0\\
\{S_1,S_1\}&=-2\{S_0,S_2\}\\
&\cdots
\end{align*}
The second equation in the list above tells that
$\text{ad}_{S_1}$ is no more a differential. On the other
hand, since $S_0$ does not depend on the antifields, at
critical points of $S_0$, i.e., at points $p$ in $W$ such
that
\[
\left.\frac{\partial S_0}{\partial x^i}\right\vert_p=0,
\qquad
\forall i,
\]
one has $\{S_0,S_2\}=0$. This means that
$\delta=\text{ad}_{S_1}$ is a differential when restricted
to the set of critical points of $S_0$. One refers to this
phenomenon by saying that $\delta$ is differential
\emph{on-shell}. When $\delta$ is a differential on the whole
of
$W$, as in the case considered before, one says that
$\delta$ is a differential \emph{off-shell}. 
\section{Traceless representations and the quantum master
equation}
In the above section we have seen how to build a solution to
the master equation from a function $S_0$
and a Lie algebra of infinitesimal symmetries for this
function. We now show how, if the representation $\rho\colon
\g\to V$ inducing the infinitesimal symmetries of the action
is traceless, then the function $S=S_0+\hbar S_1$ is actually
a solution to the quantum master equation.
\par
Recall that the quantum master equation is the equation
\[
2\U\hbar\Delta S-\{S,S\}=0
\]
If $S=S_0+ \hbar S_1$ as in the previous section,
the quantum master equation is reduced to
\[
\Delta S_1=0
\]
It is immediate to compute
\[
\Delta S_1=\Delta(x_i^+\delta
x^i)=\frac{\partial}{\partial x_i}(\delta
x^i)=\text{div}\delta
\]
i.e., $\Delta S_1$ is the divergence of the vector field
$\delta$; since $\delta$ is a vector field on $V\oplus
\Pi\g$, its divergence is a function on this space. Now,
recall that
\[
\delta=\rho^i_{jk}v^jc^k\frac{\partial}{\partial v^i}+
\frac{1}{2}f^i_{jk}c^jc^k\frac{\partial}{\partial c^i}
\]
so that
\[
\text{div}\delta=(\rho^i_{ik}+f^i_{ik})c^k,
\]
i.e.,
\[
\delta\colon (v,\g)\mapsto
\text{Tr}(\text{ad}_g)+\text{Tr}(\rho_g)
\]
Therefore, if the adjoint representation of $\g$ on itself
and the representation $\rho$ are traceless, then
$S=S_0+\hbar S_1$ is a solution to the quantum master
equation. Since $S\bigr\vert_V=S_0\bigr|_V$, we have
achieved our goal to extend $e^{\frac{i}{\hbar}S_0}$ to a
$\Delta$-closed function.

\vbox{

\medskip

\noindent{\sc Dipartimento di Matematica ``Guido Castelnuovo'' ---
Universit\`a degli Studi di Roma ``la Sapienza'' --- P.le Aldo Moro, 2 --
00185 -- Roma, Italy}
\par
\noindent{\it E-mail address:} {\tt fiorenza@mat.uniroma1.it}
}
\end{document}